\documentclass[letterpaper, 10pt, conference] {ieeeconf}
\IEEEoverridecommandlockouts                              
\usepackage[pdftex]{graphicx}
\usepackage{amsmath, amssymb}
\usepackage{placeins}
\usepackage{authblk} 
\usepackage [top=72pt,left=54pt,right=54pt,bottom=54pt] {geometry}

\usepackage{amsthm}

\newcommand {\qol} [1] {\overline{#1}}

\newcommand {\qc} [1] {\check{#1}}
\newcommand {\qR} [0] {\mathbb{R}}

\title {Globally Stable Attitude Control along Neutrally Stable Trajectories}
\author {Yujendra Bharathi Mitikiri$^1$, {\it Senior Member, IEEE}
\thanks{$^*$This work is supported by SERB Startup Research Grant SRG/2022/001423}
\thanks{$^1$ The authors are with the Department of Mechanical Engineering, Indian Institute of Technology, Tirupati, India. Email: yujendra@iittp.ac.in}
}

\begin{document}

\pagenumbering {gobble}

\maketitle

\begin{abstract}
It is quite often claimed, and correctly so, that linear methods cannot achieve global stability results for attitude control, and conversely that nonlinear control is essential in order to achieve (almost) globally stable tracking of general attitude trajectories. On account of this definitive result, and also because of the existence of powerful nonlinear control techniques, there has been relatively very little work analyzing the limits and performance of linear attitude control. This paper provides a characterization of the stability achievable for one class of linear attitude control problems, namely those leading to a constant quaternion difference, and shows how the strategy can be used to design an efficient attitude controller along a neutrally stable trajectory. In this paper, we analytically derive a critical error angle below which linearized dynamics lead to natural neutral stability for such a system, and above which the system is unstable. The dynamics are then used to derive a locally stable linear attitude controller and a globally stable nonlinear attitude controller whose performance is validated using simulations. The globally stable nonlinear controller is also shown to provide more efficient control compared to previously published work. 
\end{abstract}

\date{2026 March 11}

\section {Introduction} \label{sec:intro}

\begingroup
\allowdisplaybreaks

The attitude of a vehicle may be controlled to track a desired trajectory beginning from an (almost) arbitrary initial condition using a variety of (almost) globally stable nonlinear control algorithms, such as \cite{wen:91a}--\cite{shi:23a}. The controller in \cite{wen:91a} uses a quaternion error variable to derive the controller, and this is extended in several works such as \cite{tsiotras:95a} to axisymmetric bodies, in \cite{fortuna:01a} for neural network applications, and in \cite{mitikiri:22a} to account for modelled dissipation. An orthogonal matrix based representation of the attitude error is used in \cite{chaturvedi:11a} to achieve similar stability results, and has been extended in several works such as \cite{lee:13a} for robust control, \cite{zlotnik:14a} for output-feedback attitude control, \cite{chen:19a} for flexibile bodies, and \cite{shi:23a} for constrained systems.

Despite the known limitations of linear attitude control, its use is still widespread on account of the associated simplicity, and linear response to principal angle errors. For example, \cite{dwyer:84a} used a linearization of the double integral form of attitude dynamics to regulate large initial errors. A controller that is linear in terms of the quaternion vector component is described in \cite{bach:93a}, which has the associated limitation of windup when the quaternion error angle is greater than 180 degrees. Linear PD and LQR methods for controlling a spacecraft's Euler-angle attitude are compared in \cite{beatty:06a} leading to the conclusion that LQR methods lower the required control effort for similar attitude errors. In \cite{stingu:09a}, the authors describe a linear attitude controller using Euler angle errors for a quadrotor UAV, and the controller works as long as the quadrotor's attitude is sufficiently far away from the gimbal-lock singularity. In fact, such linear controllers are still prevalent in many commercial jet airplane autopilots which operate for the most part far away from the singular Euler angle attitude. Other similar examples from literature include a design using the truncated Taylor series approximation for the quaternion exponentials in order to achieve a fast single-axis reorientation slew-rate on agile spacecrafts in \cite{gupta:98a}, a PID attitude controller in \cite{wie:02a}, and an optimal linear controller in \cite{ghiglino:15a}.

In this paper, we are interested in analyzing the stability of linear quaternion differential equations in time. In particular we consider the special case of a constant quaternion difference, since the results for this special case are sufficiently novel (the results have not been published previously in the literature to the best of the author's knowledge) and remarkable (being one of the few problems whose solutions can be analyzed and solved by hand). The chief contributions of this article consist of a novel stability analysis of the constant quaternion difference system, and the presentation of an efficient control design based on that analysis.

\section {Notation and Terminology} \label{sec:notation}

The attitude of a rigid body will be represented by a unit-magnitude 4-component quaternion $\qc q = [q_0\quad q_1\quad q_2\quad q_3]^\text{T}$, with a check accent above the symbol representing it. The first component $q_0$ of $\qc q$ is also referred to as the scalar component $q_{\text s}$, while the three components $[q_1\quad q_2\quad q_3]^\text{T}$ are also referred to as the vector part $\vec q_\text v$. The most important law of quaternion algebra is their multiplication rule \cite{PhillipsWF:10a}:
\begin{align}
\qc p\otimes\qc q & = \begin{bmatrix} p_0 \\ \vec p_{\text v} \end{bmatrix} \otimes \begin{bmatrix} q_0 \\ \vec q_{\text v} \end{bmatrix} = \begin{bmatrix} p_0q_0 - \vec p_{\text v}^\text{T}\vec q_{\text v} \\ p_0\vec q_{\text v} + q_0\vec p_{\text v} + \vec p_{\text v}\times \vec q_{\text v} \end{bmatrix}, \label{eqn:pXq}
\end{align}
where $\vec p_{\text v}\times \vec q_{\text v}$ (third term in the vector part of the product) denotes the usual vector product of two Euclidean vectors in three-dimensional Euclidean space. The skew-symmetric matrix corresponding to the vector cross product will be denoted as $[\vec p_{\text v}\times]$ so that $\vec p_{\text v}\times \vec q_{\text v}$ may also be expressed in matrix notation as $[\vec p_{\text v}\times]\vec q_{\text v}$. It follows from the above multiplication rule that the multiplicative identity is the quaternion with unit scalar part and zero vector part, and that the inverse of a unit-magnitude quaternion is given by negating its vector part \cite{PhillipsWF:10a}:
\begin{align}
\qc q & = \begin{bmatrix} q_0 \\ \vec q_{\text v} \end{bmatrix} \Rightarrow \qc q^{-1} = \begin{bmatrix} q_0 \\ -\vec q_{\text v} \end{bmatrix}. \label{eqn:qinv}
\end{align}
We will refer to unit-magnitude quaternions simply as unit quaternions for brevity.

It may be recalled that unit quaternions provide a minimal, singularity-free, and global description of rotations \cite{Shuster:93a}, with composition of rotations achieved by multiplication of quaternions. The kinematics of an attitude quaternion may be expressed in terms of the body-frame angular velocity $\vec\omega$ as \cite{PhillipsWF:10a}
\begin{align}
\dot{\qc q} & = \frac{1}{2} \qc q\otimes\vec\omega. \label{eqn:qdot_qomg}
\end{align}
When Euclidean 3-vectors such as $\vec\omega$ are used as quaternion operands as in \eqref{eqn:qdot_qomg}, they are prefixed with a scalar component of zero.

User control over rigid-body rotation is ultimately exercised not through the angular velocity itself, but rather through torques (or moments of forces) which influence the angular acceleration and hence the angular velocity. The relationship between the angular velocity and the torque may be described by Euler's equations in vector form as
\begin{align}
\vec n & = \vec\omega\times J\vec\omega + J\dot{\vec\omega}, \label{eqn:n_Jomgdot}
\end{align}
where $\vec n \in \qR^3$ is the net external torque, and $J \in \qR^{3\times 3}$ is the (second) moment of inertia of the rigid body in the body-fixed coordinate frame. Assuming a solid body, the moment of inertia is positive definite, and the above equation may be inverted to obtain a one-to-one correspondence between the angular acceleration and the applied external torque. Therefore, one may justifiably consider the angular acceleration itself as the user input instead of the torque.

The problem we are concerned with in this article is that of controlling the attitude $\qc q$ of a rigid body so as to track a reference attitude trajectory $\qc p(t)$, whose kinematics are specified to be
\begin{align}
\dot{\qc p} & = \frac{1}{2} \qc p\otimes\vec\nu, \label{eqn:pdot_pnu}
\end{align}
where $\vec\nu$ is the reference angular velocity. We will assume that the angular acceleration $\dot{\vec\nu}$ is also given to us. We will sometimes use the half-angular-velocities, $\vec v$ for $\vec\nu/2$ and $\vec w$ for $\vec\omega/2$, to avoid the factors of 2 and reduce notational clutter. In terms of $\vec v$ and $\vec w$, we have
\begin{align}
\dot{\qc p} & = \qc p \otimes \vec v,\; \dot{\qc q} = \qc q \otimes \vec w. \label{eqn:pdot_qdot_pvqw}
\end{align}

\section{Linear constant difference attitude dynamics} \label{sec:dyn}

In this section, we will consider the linear stability of the constant difference attitude control problem. The problem is interesting on its own from a theoretical standpoint. From a practical perspective, the problem manifests, for example, in the presence of an unknown disturbance that is quasistatic (constant, or very slowly varying relative to time scales relevant to the problem) referred to the body-frame. In the presence of such a disturbance, we cannot ensure convergence of the attitude difference to zero. Instead, it is its derivative alone which can be driven to zero. 

\subsection{Nominal Dynamics}

From a purely geometric perspective, the constant difference constraint implies that the angle between the two attitudes is also a constant. Let
\begin{align}
\qc r & = \qc q - \qc p \label{eqn:r_qp}
\end{align}
denote the constant difference between the desired attitude $\qc p$ and the body attitude $\qc q$.
Define the unit error quaternion as
\begin{align}
\qc e & = \qc p^{-1}\otimes \qc q. \label{eqn:e_pinvq}
\end{align}

The scalar part of the error rotation may be expressed in terms of $\qc r$ and $\qc p$ as
\begin{align}
e_0 & = p_0q_0 + \vec p_{\text v}^\text{T}\vec q_{\text v} = p_0(p_0 + r_0) + \vec p_{\text v}^\text{T}(\vec p_{\text v} + \vec r_{\text v}) \nonumber\\
 & = 1 + r_0p_0 + \vec r_{\text v}^\text{T}\vec p_{\text v}. \label{eqn:e0_pinvr0}
\end{align}
Similarly, $e_0$ may also be expressed in terms of $\qc r$ and $\qc q$ as
\begin{align}
e_0 & = 1 - r_0q_0 - \vec r_{\text v}^\text{T}\vec q_{\text v}. \label{eqn:e0_qinvr0}
\end{align}
By adding up \eqref{eqn:e0_pinvr0} and \eqref{eqn:e0_qinvr0}, and dividing by two, we obtain
\begin{align}
e_0 & = 1 - |\qc r|^2/2. \label{eqn:e0_rmagsqrby2}
\end{align}
The constant difference constraint \eqref{eqn:r_qp} may now be invoked to conclude that $e_0$, and therefore the error angle, are also constant in time.

The difference quaternion $\qc r$ is clearly not a unit quaternion in general, but that does not preclude linear control based upon it. We are therefore interested in studying the linear (small perturbation) stability of solutions to the quaternion differential equation
\begin{align}
\dot{\qc q} = \dot{\qc p} . \label{eqn:qdot_pdot}
\end{align}

Using \eqref{eqn:pdot_qdot_pvqw} in \eqref{eqn:qdot_pdot}, we obtain the mathematical statement which encapsulates the constant difference constraint as a relationship between the desired and body (half) angular velocities, $\vec v$ and $\vec w$:
\begin{align}
\qc p\otimes\vec v = \qc q\otimes\vec w \; \Leftrightarrow \; \vec v & = \qc e \otimes \vec w, \nonumber 
\end{align}
or in terms of the scalar and vector parts,
\begin{align}
\vec e_{\text v}^\text{T}\vec w & = 0, \label{eqn:evTw_0}\\
\vec v & = e_0\vec w + \vec e_{\text v}\times \vec w. \label{eqn:v_e0wevxw}
\end{align}
It may be noted that we could have as well used the constant difference constraint and the inverse of \eqref{eqn:e_pinvq}, to infer that $\vec v$ is also perpendicular to $\vec e_{\text v}$ like $\vec w$, and to express $\vec w$ in terms of $\vec v$ as follows
\begin{align}
\vec w & = \qc e^{-1} \otimes \vec v \nonumber\\
& \Rightarrow \vec e_{\text v}^\text{T}\vec v = 0, \label{eqn:evTv_0}\\
& \hphantom{\Rightarrow} \vec w = e_0\vec v - \vec e_{\text v}\times \vec v. \label{eqn:w_e0vevxv}
\end{align}
The former set of equations, \eqref{eqn:evTw_0} and \eqref{eqn:v_e0wevxw}, will however turn out to be of greater utility for us in the following stability analysis.

The relationship in \eqref{eqn:v_e0wevxw} may be utilized to eliminate the desired (half) angular velocity $\vec v$ in the error dynamics, which may be derived using \eqref{eqn:e_pinvq} and \eqref{eqn:pdot_qdot_pvqw} as follows:
\begin{align}
\dot{\qc e} & = -\vec v \otimes \qc p^{-1} \otimes \qc q + \qc p^{-1} \otimes \qc q \otimes \vec w \nonumber\\
 & = \qc e \otimes \vec w - \vec v \otimes \qc e. \label{eqn:edot_weev}
\end{align}
Now, using \eqref{eqn:v_e0wevxw} to substitute for $\vec v$, we may express the kinematics of the error quaternion solely in terms of its own components $e_0$ and $\vec e_{\text v}$ and the body (half) angular velocity $\vec w$:
\begin{align}
\dot{\qc e} & = \begin{bmatrix} \vec e_{\text v}^\text{T}(\vec v - \vec w) \\ e_0(\vec w - \vec v) + \vec e_{\text v}\times (\vec v + \vec w) \end{bmatrix} \label{eqn:edot_1}\\
 & = \begin{bmatrix} (e_0 - 1)\vec e_{\text v}^\text{T}\vec w \\ (e_0 - 1)\vec w + \vec e_{\text v}\times \vec w + \vec e_{\text v}\vec e_{\text v}^\text{T}\vec w \end{bmatrix}. \label{eqn:edot_2}
\end{align}

We have so far utilized the relation in \eqref{eqn:v_e0wevxw} to derive the kinematic governing equation \eqref{eqn:edot_2}. The scalar part of \eqref{eqn:edot_2} is not independent of the vector part since $\qc e$ is constrained to have unit magnitude, and may be ignored. 
The vector part of \eqref{eqn:edot_2} effectively contains three scalar equations in terms of six scalar components, three each in the vectors $\vec e_{\text v}$ and $\vec w$. In order to complete the dynamic specification of the system, we need to bring in attitude kinetics and derive an equation for the body (half) angular acceleration $\dot{\vec w}$ which invokes the control torques and moments through \eqref{eqn:n_Jomgdot}.

Unlike the dynamics of the error quaternion $\qc e$, the scalar and vector parts of the dynamics of the difference quaternion $\qc r$ are independent, since its components are not constrained to yield unit magnitude. Therefore, both \eqref{eqn:evTw_0} and \eqref{eqn:v_e0wevxw} need to be enforced in the nominal motion. We have already exhausted the three scalar degrees of constraint contained in \eqref{eqn:v_e0wevxw} in deriving the three independent scalar equations comprising the vector part of \eqref{eqn:edot_2}. The constraint in \eqref{eqn:evTw_0} has to be incorporated yet, and to this end, we take its derivative to invoke the moment equation, and using \eqref{eqn:edot_2}, we obtain one scalar equation:
\begin{align}
\vec e_{\text v}^\text{T}\dot{\vec w} = -\vec w^\text{T}\dot{\vec e}_{\text v} = (1-e_0)|\vec w|^2 - (\vec e_{\text v}^\text{T}\vec w)^2. \label{eqn:evTwdot}
\end{align}

We may now design a nominal $\dot{\vec w}$ of minimal magnitude:
\begin{align}
\dot{\vec w} & = \frac{(1 - e_0)|\vec w|^2\vec e_{\text v}}{|\vec e_{\text v}|^2} = \frac{|\vec w|^2\vec e_{\text v}}{1 + e_0}. \label{eqn:wdot}
\end{align}
It may be noted that \eqref{eqn:wdot} reflects only the nominal dynamics of the angular velocity. In practice, one usually adds a feedback controller to asymptotically drive the system towards satisfying \eqref{eqn:evTw_0}. See Section \ref{sec:ctrldsgn} for one possible linear design that yields asymptotic convergence on the angular velocity.

\subsection{Linearized (Perturbation) Dynamics}

Dynamical equations \eqref{eqn:edot_2} and \eqref{eqn:wdot}, together with the constant angle constraint $\dot e_0 = 0$, define the nominal motion for $\qc e$ and $\vec\omega$ in the following perturbation and stability analysis. They may be written in (nonlinear) state-space form as
\begin{align}
\begin{bmatrix} \dot{\vec e}_{\text v} \\ \dot{\vec w} \end{bmatrix} & = \begin{bmatrix} (e_0 - 1)\vec w + \vec e_{\text v}\vec e_{\text v}^\text{T}\vec w + \vec e_{\text v}\times \vec w \\ |\vec w|^2\vec e_{\text v}/(1 + e_0) \end{bmatrix}. \label{eqn:edotwdot}
\end{align}
We may linearize the above dynamics, and use the constant difference constraint $\vec e_{\text v}^\text{T}\vec w = 0$, to obtain
\begin{align}
\begin{bmatrix} \dot{\tilde e}_{\text v} \\ \dot{\tilde w} \end{bmatrix} 
 & = A\begin{bmatrix} \tilde e_{\text v} \\ \tilde w \end{bmatrix}, \label{eqn:xdot_Ax}
\end{align}
where $\tilde e_{\text v}$ and $\tilde w$ are perturbations in the vector part $\vec e_{\text v}$ of the attitude error, and the (half) angular velocity $\vec w$, (so, the corresponding total quantities are $\vec e_{\text v} + \tilde e_{\text v}$ and $\vec w + \tilde w$) and

\begin{align}
A & = \begin{bmatrix} \vec e_{\text v}\vec w^\text{T} - [\vec w\times] - \dfrac{\vec w\vec e_{\text v}^\text{T}}{e_0} & (e_0 - 1)1_{3\times 3} + [\vec e_{\text v}\times] + \vec e_{\text v}\vec e_{\text v}^\text{T} \\ \dfrac{|\vec w|^21_{3\times 3}}{1 + e_0} + \dfrac{|\vec w|^2e_{\text v}\vec e_{\text v}^\text{T}}{(1 + e_0)^2e_0} & \dfrac{2\vec e_{\text v}\vec w^\text{T}}{1 + e_0} \end{bmatrix}. \label{eqn:A_1}
\end{align}
We have substituted $\tilde e_0$ with $-\vec e_{\text v}^\text{T}\tilde e_{\text v}/e_0$ in the above expression. The trace of $A$ above may be immediately verified to be equal to zero after using \eqref{eqn:evTw_0}. This implies that all the eigenvalues are on the imaginary axis, or at least one of them has a positive real part and the system is not asymptotically stable.

Although the $6\times 6$ matrix $A$ in \eqref{eqn:A_1} looks quite formidable, it may be simplified significantly by applying the following orthogonal and scaling transformations. Let
\begin{align}
C & = \begin{bmatrix} \dfrac{\vec e_{\text v}}{|\vec e_{\text v}|} & \dfrac{\vec w\times \vec e_{\text v}}{|\vec e_{\text v}||\vec w|} & \dfrac{\vec w}{|\vec w|} \end{bmatrix}; \label{eqn:C}
\end{align}
then, it is a simple matter to verify that
\begin{equation}
\begin{aligned}
C^\text{T}C & = 1_{3\times 3}, &
C^\text{T}\vec e_{\text v}\vec e_{\text v}^\text{T}C & = \begin{bmatrix} |\vec e_{\text v}|^2 & & \\ & 0 & \\ & & 0 \end{bmatrix}, \\ 
C^\text{T}\vec e_{\text v}\vec w^\text{T}C & = \begin{bmatrix} & & |\vec e_{\text v}||\vec w| \\ & 0 & \\ 0 & & \end{bmatrix}, & 
C^\text{T}\vec w\vec e_{\text v}^\text{T}C & = \begin{bmatrix} & & 0 \\ & 0 & \\ |\vec e_{\text v}||\vec w| & & \end{bmatrix}, \\ 
C^\text{T}[\vec e_{\text v}\times]C & = \begin{bmatrix} 0 \\ & & -|\vec e_{\text v}| \\ & |\vec e_{\text v}| & \end{bmatrix}, & 
C^\text{T}[\vec w\times]C & = \begin{bmatrix} & -|\vec w| \\ |\vec w| & & \\ & & 0 \end{bmatrix}. 
\end{aligned} \label{eqn:CTC_1}
\end{equation}
We now rotate the state matrix $A$ in \eqref{eqn:A_1} through $C$
\begin{align}
& A' = \begin{bmatrix} C^\text{T} \\ & C^\text{T} \end{bmatrix} A \begin{bmatrix} C \\ & C \end{bmatrix} \label{eqn:CTAC}\\
& = \left[ \begin{array} {c|c}
\begin{matrix} \hphantom{\dfrac{w^2}{c(1\!+\!c)}} & w & sw \\ -w & \hphantom{\dfrac{w^2}{1\!+\!c}} & \hphantom{\dfrac{w^2}{1\!+\!c}} \\ \!\!-sw/c & \hphantom{\dfrac{w^2}{1\!+\!c}} \end{matrix}
  & \begin{matrix} c\!-\!c^2 & & \\ & c\!-\!\!1 & -s \\ & s & c\!-\!\!1 \end{matrix}\\
\hline
\begin{matrix} \dfrac{w^2}{c(1\!+\!c)} \\ & \dfrac{w^2}{1\!+\!c} \\ & & \dfrac{w^2}{1\!+\!c} \end{matrix}
  & \begin{matrix} \hphantom{c\!-\!c^2} & & \dfrac{2sw}{1\!+\!c} \\ & \hphantom{c}0\hphantom{1}\vphantom{\dfrac{w^2}{c}} \\ 0 & & \vphantom{\dfrac{w^2}{1\!+\!c}} \end{matrix} 
\end{array} \right], \label{eqn:A_2}
\end{align}
where, $c = e_0$ and $s = |\vec e_{\text v}|$ are the cosine and sine of the error half angle, and $w = |\vec w|$. Note that empty entries in a matrix indicate a zero.

Since $A$ and $A'$ are related through the orthogonal transformation in \eqref{eqn:CTAC}, they share the same characteristic polynomial. The characteristic polynomial $\det(\lambda 1_{6\times 6} - A')$ for the matrix $A'$ in \eqref{eqn:A_2} may be evaluated to be equal to (see the appendix \ref{sec:cpderiv} for a derivation)
\begin{align}
& f_A(\lambda) = f_{A'}(\lambda) = \frac{\lambda^6}{|\vec w|^6} + \frac{(1 + 3e_0 - e_0^2 - e_0^3)}{e_0(1 + e_0)}\frac{\lambda^4}{|\vec w|^4} \nonumber\\
& \quad + \frac{(1 - e_0)(3 + 3e_0 - 2e_0^2)}{(1 + e_0)^2}\frac{\lambda^2}{|\vec w|^2} + \frac{2e_0(1 - e_0)^2}{(1 + e_0)^3}. \label{eqn:lamcharpoly}
\end{align}
An immediate observation is that the characteristic polynomial $f_A(\lambda)$ contains only even powers of $\lambda$ so roots are point symmetric about the origin. Multiplying $f_A$ with $(1 + e_0)^3$ and substituting $\lambda'$ for $(1 + e_0)\lambda^2/|\vec w|^2$, $f_A(\lambda)$ may be simplified to the following cubic in terms of $\lambda'$
\begin{align}
& \lambda'^3 + (1/e_0 + 3 - e_0 - e_0^2)\lambda'^2 \nonumber\\
& \quad + (1 - e_0)(3 + 3e_0 - 2e_0^2)\lambda' + 2e_0(1 - e_0)^2. \label{eqn:lamprimecharpoly}
\end{align}
Since the roots of a cubic polynomial can be obtained in closed form, it is therefore possible to analytically compute the eigenvalues of the sixth order system in \eqref{eqn:xdot_Ax}. This rather remarkable outcome is a consequence of the special structure associated with the dynamical equations \eqref{eqn:edot_2} and \eqref{eqn:wdot} of the system under consideration \eqref{eqn:xdot_Ax}, and their linearization about the constraint $\vec e_{\text v}^\text{T}\vec w = 0$ in \eqref{eqn:evTw_0}. In particular, the elements of $A$ may all be expressed in terms of standard Euclidean vector operations such as the scalar and vector product, and the projection operator. And the characteristics of systems described using such Euclidean vector expressions are invariant upon the use of orthogonal transformations. The specific orthogonal transformation $C$ used in \eqref{eqn:CTAC} aligns the $x$-axis along $\vec e_{\text v}$ and the $z$-axis along $\vec w$.

The roots of the cubic in $\lambda'$ are all real numbers when the non-trivial factor of the discriminant
\begin{align}
\Delta & = 1 + 5e_0 - 8e_0^2 + 4e_0^3 \label{eqn:Delta}
\end{align}
is positive (the trivial factors consist of a root at $-1$ and eight at $+1$). Within the range $[-1,\; 1]$ for $e_0$, the discriminant is positive when $e_0 \gtrsim -0.16$ (error angle less than $\approx 99^\circ$) and negative otherwise. Furthermore, the coefficients in $f_A$ are nonnegative for all positive values of $e_0$, and the product of the coefficients of $\lambda'^2$ and $\lambda'$ is greater than the constant coefficient in \eqref{eqn:lamprimecharpoly}, allowing us to conclude that the roots of \eqref{eqn:lamprimecharpoly} are always in the left half of the complex plane (that is, their real parts are negative). When they are purely real negative, that leads to purely imaginary roots in \eqref{eqn:lamcharpoly}, and the system is neutrally stable despite the lack of active negative feedback induced stabilization. When the roots of \eqref{eqn:lamprimecharpoly} are complex for $e_0 \lesssim -0.16$ (error angle $\gtrsim 99^\circ$), that leads to \eqref{eqn:lamcharpoly} having roots with both positive and negative real parts and consequent instability. 

The presence of purely imaginary eigenvalues for the system in \eqref{eqn:xdot_Ax} for bounded error angles tells us that the system is neutrally stable and does not diverge away from its nominal difference, as long as the constraints in \eqref{eqn:evTw_0} and \eqref{eqn:v_e0wevxw} are satisfied. This justifies designing a continuous switching controller comprising of a globally stable nonlinear controller for error angles above a certain bound, and a locally stable linear controller for small errors. 

\section {Attitude Control Design} \label{sec:ctrldsgn}

\subsection {Nominal Trajectory Parameterization} \label{sec:nomtrajparam}

Without loss of generality, a nominal constant quaternion difference trajectory for the desired and body attitudes $\qc p$ and $\qc q$ may be parameterized (modulo a constant rotation isomorphism) in terms of two scalar signals $\alpha(t)$ and $\beta(t)$, and a scalar constant $\phi$ as follows:
\begin{align}
\qc p & = \begin{bmatrix} c_\phi c_\alpha \\ c_\phi s_\alpha c_\beta \\ c_\phi s_\alpha s_\beta \\ -s_\phi \end{bmatrix} ,\; \qc q = \begin{bmatrix} c_\phi c_\alpha \\ c_\phi s_\alpha c_\beta \\ c_\phi s_\alpha s_\beta \\ s_\phi \end{bmatrix}, \label{eqn:nompq}\\
\Rightarrow \qc e & = \qc p^{-1} \otimes \qc q \nonumber\\
 & = \begin{bmatrix} c_{2\phi} \\ -s_{2\phi} s_\alpha s_\beta \\ s_{2\phi} s_\alpha c_\beta \\ s_{2\phi} c_\alpha \end{bmatrix} = c_{2\phi} \qc 1 + s_{2\phi} \begin{bmatrix} 0 \\ -s_\alpha s_\beta \\ s_\alpha c_\beta \\ c_\alpha \end{bmatrix} . \label{eqn:nome}
\end{align}
Let us assume that $\alpha(t)$ and $\beta(t)$ are arbitrary twice-differentiable functions of time. Taking the derivatives, we obtain the desired and body (half-)angular velocities $v$ and $w$ in terms of the constant parameter $\phi$, and the signals $\alpha(t)$ and $\beta(t)$:
\begin{align}
\dot{\qc p} & = \dot{\qc q} = \dot\alpha c_\phi \begin{bmatrix} -s_\alpha \\ c_\alpha c_\beta \\ c_\alpha s_\beta \\ 0 \end{bmatrix} + \dot\beta c_\phi s_\alpha \begin{bmatrix} 0 \\ -s_\beta \\ c_\beta \\ 0 \end{bmatrix} \label{eqn:nompdotqdot}\\
\vec v & = \qc p^{-1} \otimes \dot{\qc p} \nonumber\\
 & = (\dot\alpha c_\phi^2 - \dot\beta c_\phi s_\phi s_\alpha) \begin{bmatrix} c_\beta \\ s_\beta \\ 0 \end{bmatrix} + (\dot\alpha c_\phi s_\phi + \dot\beta c_\phi^2 s_\alpha) \begin{bmatrix} -c_\alpha s_\beta \\ c_\alpha c_\beta \\ -s_\alpha \end{bmatrix} . \label{eqn:nomv}
\end{align}
Since $\qc q$ and $\vec w$ are obtained from $\qc p$ and $\vec v$ by replacing $\phi$ with its negative $-\phi$, we must have
\begin{align}
\vec w & = (\dot\alpha c_\phi^2 + \dot\beta c_\phi s_\phi s_\alpha) \begin{bmatrix} c_\beta \\ s_\beta \\ 0 \end{bmatrix} + (\dot\beta c_\phi^2 s_\alpha - \dot\alpha c_\phi s_\phi) \begin{bmatrix} -c_\alpha s_\beta \\ c_\alpha c_\beta \\ -s_\alpha \end{bmatrix}. \label{eqn:nomw}
\end{align}
We can see from \eqref{eqn:nome}, \eqref{eqn:nomv} and \eqref{eqn:nomw} that $\vec e_{\text v}^\text{T}\vec v = \vec e_{\text v}^\text{T}\vec w = 0$, which motivates us to define the orthonormal triad
\begin{align}
\begin{bmatrix} \hat x & \hat y & \hat z \end{bmatrix} & = \begin{bmatrix} c_\beta & -c_\alpha s_\beta & -s_\alpha s_\beta \\ s_\beta & c_\alpha c_\beta & s_\alpha c_\beta \\ 0 & -s_\alpha & c_\alpha \end{bmatrix} \nonumber\\
 & = \begin{bmatrix} c_\beta & -s_\beta \\ s_\beta & c_\beta \\ & & 1 \end{bmatrix} \begin{bmatrix} 1 \\ & c_\alpha & s_\alpha \\ & -s_\alpha & c_\alpha \end{bmatrix} , \label{eqn:nomxyzhat}
\end{align}
with derivatives
\begin{align}
\dot{\hat x} & = \dot\beta(c_\alpha \hat y + s_\alpha \hat z), \label{eqn:xhatdot}\\
\dot{\hat y} & = -\dot\alpha\hat z - \dot\beta c_\alpha \hat x, \label{eqn:yhatdot}\\
\dot{\hat z} & = -\dot\beta s_\alpha \hat x + \dot\alpha\hat y, \label{eqn:zhatdot}\\
\vec\omega & = -\dot\alpha \hat x - \dot\beta s_\alpha \hat y + \dot\beta c_\alpha \hat z. \label{eqn:nomomg}
\end{align}
From \eqref{eqn:nomv} and \eqref{eqn:nomw}, the magnitudes of the nominal $\vec v$ and $\vec w$ are
\begin{align}
\vec v^\text{T}\vec v & = \vec w^\text{T}\vec w = (\dot\alpha c_\phi^2 \mp \dot\beta c_\phi s_\phi s_\alpha)^2 + (\dot\beta c_\phi^2 s_\alpha \pm \dot\alpha c_\phi s_\phi)^2 \nonumber\\
 & = \dot\alpha^2 c_\phi^2 + \dot\beta^2c_\phi^2s_\alpha^2. \label{eqn:nomvwmag}
\end{align}
The desired (half-) angular acceleration $\dot{\vec v}$ may also be expressed in terms of the parameter $\phi$ and the signals $\alpha(t)$ and $\beta(t)$ as
\begin{align}
\dot{\vec v} & = (\ddot\alpha c_\phi^2 - 2\dot\alpha \dot\beta c_\phi s_\phi c_\alpha - \ddot\beta c_\phi s_\phi s_\alpha - \dot\beta^2 c_\phi^2 c_\alpha s_\alpha)\hat x \nonumber\\
 & \; + (\ddot\alpha c_\phi s_\phi + 2\dot\alpha \dot\beta c_\phi^2 c_\alpha + \ddot\beta c_\phi^2 s_\alpha - \dot\beta^2 c_\phi s_\phi c_\alpha s_\alpha)\hat y \nonumber\\
 & \; - (\dot\alpha^2 + \dot\beta^2 s_\alpha^2) c_\phi s_\phi \hat z. \label{eqn:nomvdot}
\end{align}
The body (half-) angular acceleration $\dot{\vec w}$ may be obtained by substituting $\phi$ with $-\phi$ in \eqref{eqn:nomvdot}.

\subsection {Feedback and Control of General Attitude Trajectories} \label{sec:fdbkctrl}

For general desired attitude trajectories that do not satisfy the requirement in \eqref{eqn:evTv_0} for constant quaternion difference dynamics, we cannot directly use feedforward parameterization and control of the body attitude presented in the previous subsection. However, the conditions in \eqref{eqn:evTw_0} and \eqref{eqn:w_e0vevxv} can still be used advantageously to design an efficient attitude controller along a neutrally stable trajectory as described in this subsection.

Since equations \eqref{eqn:evTw_0} and \eqref{eqn:w_e0vevxv} together impose four scalar degrees of constraint on the three scalar coordinates of $w$, we could use feedback to define two independent control laws for $\dot{\vec w}$
\begin{align}
\vec e_{\text v}^\text{T} \dot{\vec w} & = g(\vec v, \dot{\vec v}, \qc e, \vec w) \nonumber\\
 & = -\vec w^\text{T}\dot{\vec e}_{\text v} - k\vec e_{\text v}^\text{T}\vec w \nonumber\\
 & = e_0\vec w^\text{T}(\vec v - \vec w) - \vec e_{\text v}^\text{T} \vec v \times \vec w  - k\vec e_{\text v}^\text{T}\vec w, \label{eqn:ctrl_evTwdot}\\
\dot{\vec w} & = \vec h(\vec v, \dot{\vec v}, \qc e, \vec w) \nonumber\\
 & = \dot e_0 \vec v + e_0\dot{\vec v} - \dot{\vec e}_{\text v} \times \vec v - \vec e_{\text v} \times \dot{\vec v} - L(\vec w - e_0\vec v + \vec e_{\text v} \times \vec v) \nonumber\\
 & = e_0\dot{\vec v} - \vec e_{\text v} \times \dot{\vec v} - \vec v\, \vec e_{\text v}^\text{T} \vec w - \vec w\, \vec e_{\text v}^\text{T} \vec v + e_0\vec v\times\vec w + \vec e_{\text v} \vec v^\text{T} (\vec v + \vec w) \nonumber\\
 & \quad - L(\vec w - e_0\vec v + \vec e_{\text v} \times \vec v). \label{eqn:ctrl_wdot}
\end{align}
where $k \in \qR$ is a positive feedback gain and $L \in \qR^{3\times 3}$ is a positive definite feedback matrix. We may design the overall control input as follows: 
\begin{align}
\dot{\vec w}^\star & = \vec h + u\frac{(g - \vec e_{\text v}^\text{T}\vec h)\vec e_{\text v}}{\vec e_{\text v}^\text{T}\vec e_{\text v}}, \label{eqn:ctrl_wastdot}
\end{align}
where $u$ is a tunable weight. When $u = 0$, we identically have $\vec e_{\text v}^\text{T}\dot{\vec w}^\ast = 0$, and when $u = 1$, we identically have $\dot{\vec w}^\ast = \vec h$.
One remaining flaw with the control law in \eqref{eqn:ctrl_wastdot} is that $g$ (defined in \eqref{eqn:ctrl_evTwdot}) could remain nonzero as $\qc e$ approaches identity and $\vec e_{\text v}$ approaches zero, leading to unbounded control for $\dot{\vec w}$ in \eqref{eqn:ctrl_wastdot}. A simple remedy is to replace $\vec w$ with $e_0\vec v - \vec e_{\text v}\times\vec v$ in $\dot{\vec e}_{\text v}$ in $g$ as follows
\begin{align}
g & = e_0\vec w^\text{T}((1 - e_0)\vec v + \vec e_{\text v}\times \vec v) - \vec e_{\text v}^\text{T} \vec v \times \vec w  - k\vec e_{\text v}^\text{T}\vec w, \label{eqn:ctrl_evTwdot_2}
\end{align}
yielding a bounded control law even when $\qc e \rightarrow \qc 1$.

Rather than take the weighted average of two acceleration control laws as done above, a theoretically more fruitful approach yielding a strong stability result is to first weigh the constraints in \eqref{eqn:evTw_0} and \eqref{eqn:w_e0vevxv}, and then differentiate the combination to derive the acceleration control law as done below.

Let the target angular velocity be
\begin{align}
\vec w^\ast & = e_0\vec v - \vec e_{\text v} \times \vec v - u \frac{e_0\vec e_{\text v}\vec e_{\text v}^\text{T}\vec v} {|\vec e_{\text v}|} \label{eqn:ctrl_wast}
\end{align}
where $u$ is a weighting factor between the constraints in \eqref{eqn:evTw_0} and \eqref{eqn:w_e0vevxv}. Note that \eqref{eqn:evTw_0} is only approximately satisfied for $|\vec e_{\text v}| \neq 1$. Equation \eqref{eqn:evTw_0} could have been identically satisfied by using the combination
\begin{align}
\vec w^\ast & = e_0\vec v - \vec e_{\text v} \times \vec v - u \frac{e_0\vec e_{\text v}\vec e_{\text v}^\text{T}\vec v} {|\vec e_{\text v}|^2} \label{eqn:ctrl_wast_2}
\end{align}
with $|\vec e_{\text v}|^2$ replacing the $|\vec e_{\text v}|$ in the denominator of the third term on the right hand side of \eqref{eqn:ctrl_wast}, but this new law has the even bigger problem of diverging when $|\vec e_{\text v}| \rightarrow 0$ upon differentiation, which is an important consideration for nonlinear control design.

Returning to \eqref{eqn:ctrl_wast}, we see that
\begin{align}
\lim_{u\rightarrow 0} \vec e_{\text v}^\text{T}\vec w^\ast & = e_0(1 - |\vec e_{\text v}|) \vec e_{\text v}^\text{T} \vec v, \nonumber\\
\lim_{u\rightarrow 1} \vec w^\ast & = e_0\vec v - \vec e_{\text v} \times \vec v, \nonumber
\end{align}
which satisfy \eqref{eqn:evTw_0} approximately and \eqref{eqn:w_e0vevxv} exactly. The residual noncompliance of \eqref{eqn:evTw_0} is of the same order of magnitude as the residual noncompliance of \eqref{eqn:evTv_0} in the desired trajectory. The acceleration control law from \eqref{eqn:ctrl_wast} is therefore
\begin{align}
\dot w & = \dot e_0\vec v + e_0\dot{\vec v} - \dot{\vec e}_{\text v} \times \vec v - \vec e_{\text v} \times \dot{\vec v} - u \frac{\vec e_{\text v}\vec e_{\text v}^\text{T}} {|\vec e_{\text v}|} (\dot e_0\vec v + e_0\dot{\vec v}) \nonumber\\
 & \quad - u \frac{e_0(\dot{\vec e}_{\text v}\vec e_{\text v}^\text{T} + \vec e_{\text v}\dot{\vec e}_{\text v}^\text{T})\vec v} {|\vec e_{\text v}|} + u \frac {e_0(\vec e_{\text v}^\text{T}\dot{\vec e}_{\text v})\vec e_{\text v}\vec e_{\text v}^\text{T}\vec v} {|\vec e_{\text v}|^3} \nonumber\\
 & \quad - L \left[\vec w - e_0\vec v + \vec e_{\text v} \times \vec v + u \frac {e_0\vec e_{\text v}\vec e_{\text v}^\text{T}\vec v} {|\vec e_{\text v}|} \right] . \label{eqn:ctrl_wastdot_2}
\end{align}
In the above control law, $\dot e_0$ and $\dot{\vec e}_{\text v}$ are given in terms of $e_0$, $\vec e_{\text v}$, $\vec v$, and $\vec w$ through \eqref{eqn:edot_1}.

Calling the right hand side of \eqref{eqn:ctrl_wast} as $\vec f$, we have
\begin{align}
\dot{\vec w} & = \dot{\vec f} - L(\vec w - \vec f). \nonumber
\end{align}
Defining a positive definite candidate Lyapunov function
\begin{align}
V_L & = \frac{1}{2} |\vec w - \vec f|^2, \label{eqn:VL}
\end{align}
we obtain for its first derivative
\begin{align}
\dot V_L & = -(\vec w - \vec f)^\text{T} L (\vec w - \vec f) , \label{eqn:VLdot}
\end{align}
which is a negative definite function, and for its second derivative
\begin{align}
\ddot V_L = -(\vec w - \vec f)^\text{T}L(\dot{\vec w} - \dot{\vec f}), \nonumber
\end{align}
which is bounded and therefore implies asymptotic convergence of $\vec w$ to $\vec f$ using Barbalat's Lemma. It should now be clear why we insisted on choosing the bounded control law of \eqref{eqn:ctrl_wastdot_2}.

\section {Validation using Simulations} \label{sec:results}

In this section, we present simulation results to validate the key ideas presented in the previous section. The stability of the constant quaternion difference system is characterized by the eigenvalues of the state matrix in \eqref{eqn:A_1}, which are the roots of the characteristic polynomial in \eqref{eqn:lamcharpoly}. As we saw in the previous section, the best we can hope with respect to stability is that all the eigenvalues lie on the imaginary axis (since the trace of \eqref{eqn:A_1} is zero). In order that the $\lambda$ of \eqref{eqn:lamcharpoly} lie on the imaginary axis, the $\lambda'$ of \eqref{eqn:lamprimecharpoly} must be negative real numbers. The latter condition may be verified by evaluating the discriminant (for realness) and applying the Routh-Hurwitz test (for negativity).

The first result is to demonstrate that all roots of \eqref{eqn:lamprimecharpoly} are in the left-half of the complex plane for $e_0 > 0$. This can be easily accomplished using the Routh-Hurwitz test (see Routh-Hurwitz stability criterion in, for example, \cite{Dorf:11a}). For $0 < e_0 \le 1$, it is obvious that the coefficients of \eqref{eqn:lamprimecharpoly} are all non-negative. We also need to ensure that the product of the coefficients $a$ and $b$ of $\lambda'^2$ and $\lambda'$ is greater than the constant coefficient $c$. This condition is verified on the left side of Fig. \ref{fig:rhcoeff}. The adjacent plot on the right side shows that the discriminant $\Delta$ of \eqref{eqn:lamprimecharpoly} also remains positive for all $e_0 > 0$. Together, the plots show that we obtain neutral stability for $e_0 > 0$. For $e_0 \le 0$, the constant term in \eqref{eqn:lamprimecharpoly} is negative and we can imemdiately infer that the system is unstable.
\begin{figure} [!ht] \centering
\begin{minipage}{0.49\linewidth} \begin{center}
\includegraphics [width=\linewidth] {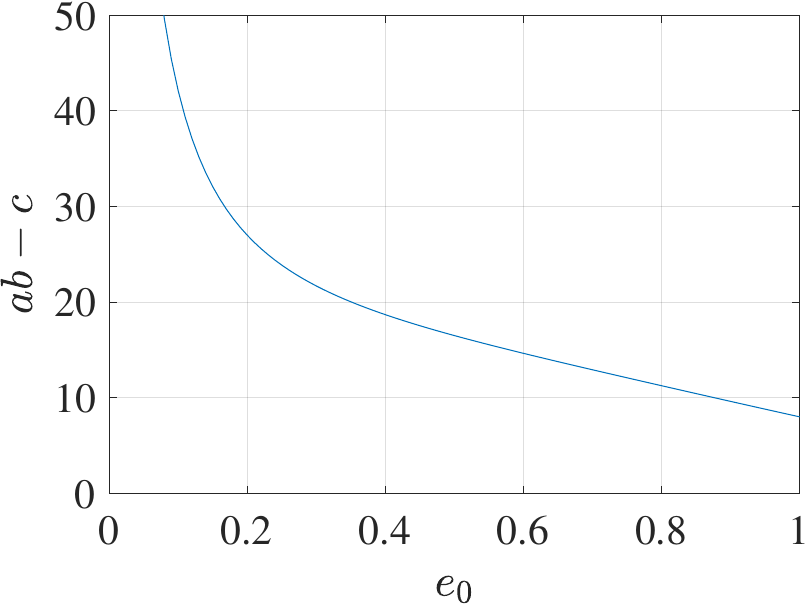}
\end{center} \end{minipage}
\begin{minipage}{0.49\linewidth}
\includegraphics [width=\linewidth] {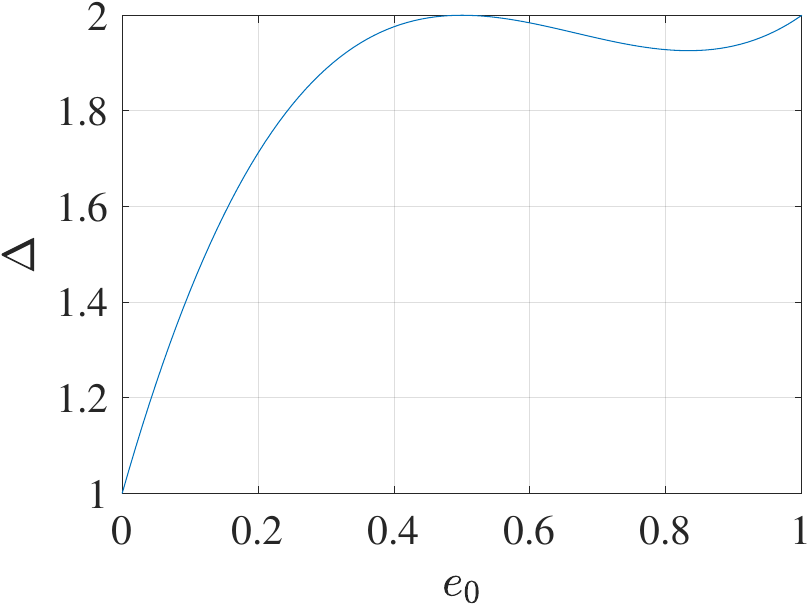}
\end{minipage}
\caption [] {{\it Left:} The product of the coefficients $a$ and $b$ of $\lambda'^2$ and $\lambda'$ in \eqref{eqn:lamprimecharpoly} must be greater than the constant coefficient $c$ in order to ensure all roots lie in the left half of the complex plane. {\it Right:} The discriminant of \eqref{eqn:lamprimecharpoly} may be plotted to graphically determine the critical error angle that separates the neutrally stable region from the unstable region for the linearized dynamics in \eqref{eqn:xdot_Ax}. Here, we see that both plots indicate neutral stability of the constant quaternion difference system for $e_0 > 0$ (error angle less than $\pi$).}
\label{fig:rhcoeff}
\end{figure}

We next verify the control design in Section \ref{sec:ctrldsgn} based upon the stability analysis in Section \ref{sec:dyn} with a specified difference of $\qc r = (0, 0, 0, 2\sin 0.2)$, which yields a nominal error $e_0 = \cos 0.4 > 0$ as required by the stability analysis. For the prescribed attitude, we use the nominal trajectory presented in Subsection \ref{sec:nomtrajparam}. For the body's attitude, we begin with an initial error of 0.6 radian with respect to the specified difference.

Below, we plot the control effort (half-angular acceleration $\dot{\vec w}$) and convergence (the Lyapunov error function $V_L$) using the controller in \cite{mitikiri:22a} or \cite{wen:91a}, \eqref{eqn:w_e0vevxv}, and \eqref{eqn:ctrl_wastdot_2}. The controller gain matrix $L$ is chosen as the identity matrix for all three controllers. The control error is plotted on the left, and the Lyapunov error on the right. Since the attitude controller using \eqref{eqn:w_e0vevxv} yields neutral stability, we see that the error remains more or less of the same magnitude (other than numerical errors) and neither converges nor diverges with time. The controller in \cite{mitikiri:22a} and \eqref{eqn:ctrl_wastdot_2} are asymptotically stable, on the other hand.
\begin{figure} [!ht] \centering
\begin{minipage}{0.49\linewidth} \begin{center}
\includegraphics [width=\linewidth] {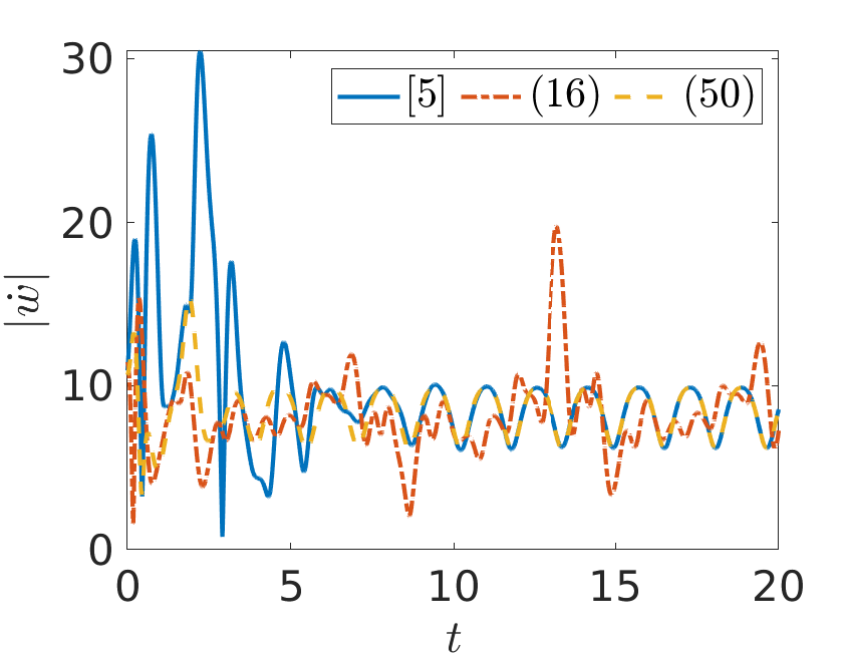}
\end{center} \end{minipage}
\begin{minipage}{0.49\linewidth}
\includegraphics [width=\linewidth] {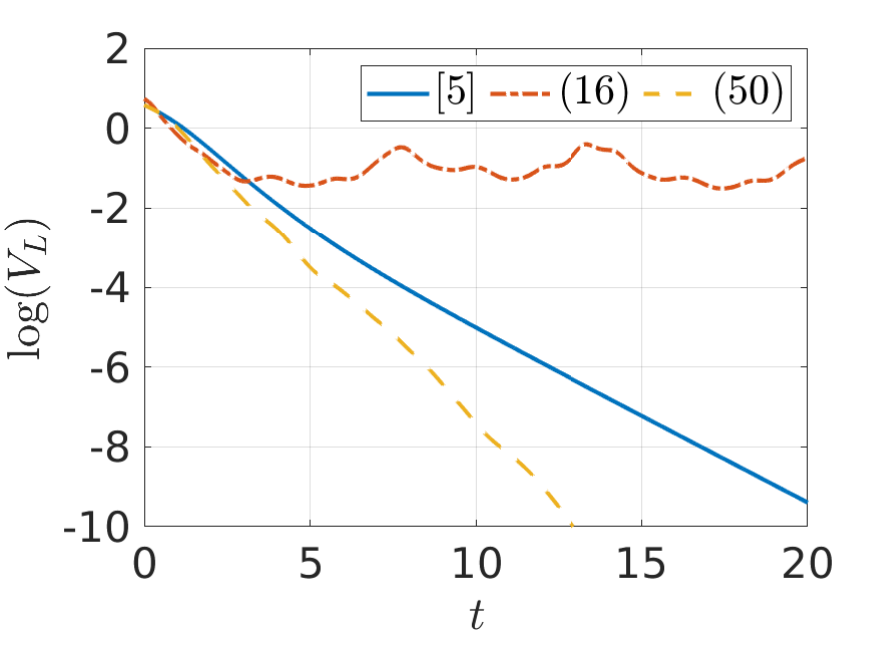}
\end{minipage}
\caption [] {{\it Left:} The controllers proposed in \eqref{eqn:w_e0vevxv} and \eqref{eqn:ctrl_wastdot_2} demand lesser effort for similar magnitude of errors when compared to the controller in \cite{mitikiri:22a}. {\it Right:} The controller proposed in \eqref{eqn:w_e0vevxv} achieves neutral stability and those proposed in \cite{mitikiri:22a} and \eqref{eqn:ctrl_wastdot_2} achieve asymptotic stability given an initial error of $0.6\pi$ rad with a nominal error $e_0 = \cos 0.4\pi$.}
\label{fig:e0pos}
\end{figure}

With a nominal error of $e_0 = \cos 0.6\pi < 0$, the behavior is very similar, except that the controller proposed in \eqref{eqn:ctrl_wastdot_2} causes the Lyapunov error function to converge at about twice as fast as that in \cite{mitikiri:22a}. This is possible because of the fact that the target angular velocity with the controller proposed in \eqref{eqn:ctrl_wastdot_2} is only half as far away from the the actual angular velocity when compared to the target angular velocity in \cite{mitikiri:22a} or \cite{wen:91a}. So an equivalent control effort (shown in the left of Fig. \ref{fig:e0neg}) leads to a stronger error convergence.
The controller in \eqref{eqn:w_e0vevxv} becomes unstable, but the error still remains bounded since the attitude space has a finite size, that is $0 \le |1 - e_0| \le 2$ always, irrespective of whether the controller is stable or not, and since \eqref{eqn:w_e0vevxv} preserves the half angular velocity magnitude.
\begin{figure} [!ht] \centering
\begin{minipage}{0.49\linewidth} \begin{center}
\includegraphics [width=\linewidth] {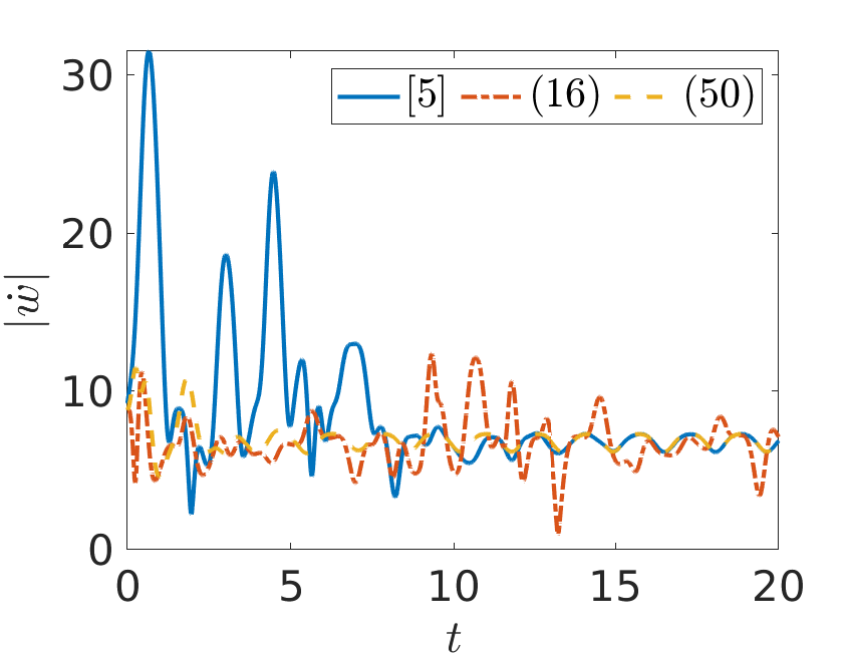}
\end{center} \end{minipage}
\begin{minipage}{0.49\linewidth}
\includegraphics [width=\linewidth] {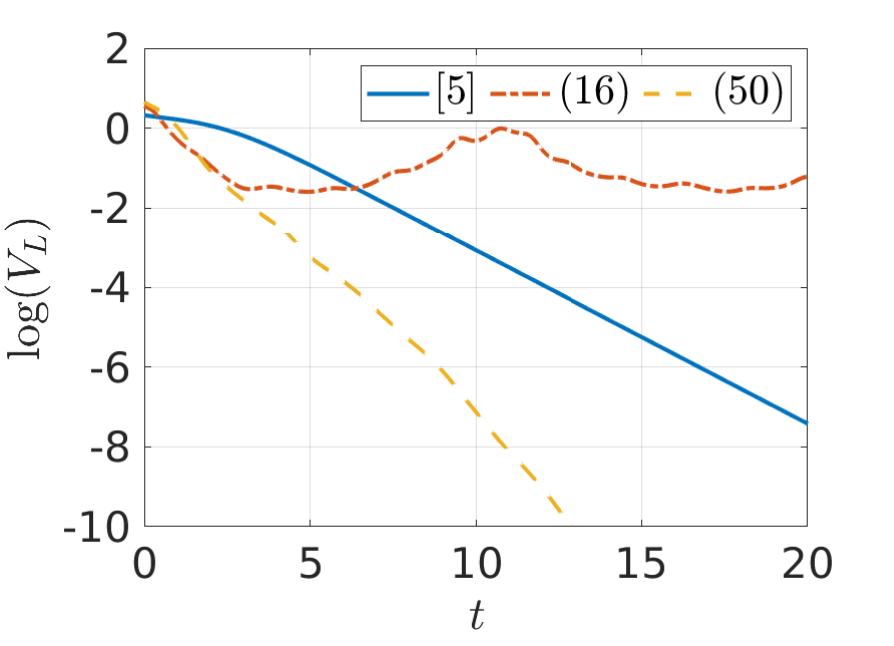}
\end{minipage}
\caption [] {{\it Left:} The controllers proposed in \eqref{eqn:w_e0vevxv} and \eqref{eqn:ctrl_wastdot_2} demand lesser effort for similar magnitude of errors when compared to the controller in \cite{mitikiri:22a}. {\it Right:} The controller proposed in \eqref{eqn:w_e0vevxv} achieves neutral stability and those proposed in \cite{mitikiri:22a} and \eqref{eqn:ctrl_wastdot_2} achieve asymptotic stability given an initial error of $0.6\pi$ rad with a nominal error $e_0 = \cos 0.6\pi$.}
\label{fig:e0neg}
\end{figure}

\section {Conclusion} \label{sec:conclude}

In this paper, we have looked at the limitations of linear attitude control for the special case of maintaining a constant quaternion difference. The resulting analysis shows that we can design a neutrally stable controller for initial error angles less than a certain value that can be derived analytically. The stability analysis can further be used to design a stable controller that achieves asymptotic convergence with lesser control effort than previous works. It would be interesting to see if it is possible to derive an optimal filtering gain $u$ in \eqref{eqn:ctrl_wastdot_2} which would yield the best compromise between the rate of asymtotic error convergence and the control cost as future work.

\section {Funding and Conflict of Interest}
This work is supported by SERB Startup Research Grant SRG/2022/001423. The authors state that there is no potential competing interest in the work presented in this article.

\bibliographystyle{unsrt}
\bibliography{ref.bib}

\appendix

\subsection {Characteristic Polynomial Derivation} \label{sec:cpderiv}

After non-dimensionalizing $A'$ in \eqref{eqn:CTAC}, we obtain the normalized matrix:
\begin{align}
& A'' = \begin{bmatrix} 1/w \\ & 1/w^2 \end{bmatrix} A' \begin{bmatrix} 1 \\ & w \end{bmatrix} \label{eqn:STAS}\\
& = \left[ \begin{array} {c|c}
\begin{matrix} \hphantom{\dfrac{1}{c(1\!+\!c)}} & 1 & s \\ -1 & \hphantom{\dfrac{1}{1\!+\!c}} & \hphantom{\dfrac{1}{1\!+\!c}} \\ \!\!-s/c & \phantom{\dfrac{1}{1\!+\!c}} \end{matrix}
  & \begin{matrix} c\!-\!c^2 & & \\ & c\!-\!\!1 & -s \\ \vphantom{\dfrac{|\vec w|}{c}} & s & c\!-\!\!1 \end{matrix}\\
\hline
\begin{matrix} \dfrac{1}{c(1\!+\!c)} \\ & \dfrac{1}{1\!+\!c} \\ & & \dfrac{1}{1\!+\!c} \end{matrix}
  & \begin{matrix} \hphantom{c\!-\!c^2} & & \dfrac{2s}{1\!+\!c} \\ & \hphantom{c}0\hphantom{1}\vphantom{\dfrac{1}{c}} \\ 0 & & \vphantom{\dfrac{1}{1\!+\!c}} \end{matrix} 
\end{array} \right]. \label{eqn:A_3}
\end{align}
Circularly rotating the top row and left column, we obtain for the determinant $|\lambda 1_{6\times 6} - A''|$,
\begin{align}
& \left| \begin{array} {c|c}
\begin{matrix} \lambda \\ & \lambda \\ & & \lambda\vphantom{-\dfrac{1}{c(1+c)}} \\ \hline \dfrac{-1}{1\!+\!c} \\ & \dfrac{-1}{1\!+\!c} \\ -1 & -s & c^2\!-\!c \end{matrix}
  & \begin{matrix} 1\!-\!c & s & 1 \\ -s & 1\!-\!c & s/c \\ & \dfrac{-2s}{1+c} & \dfrac{-1}{c(1+c)} \\ \hline \lambda \vphantom{\dfrac{-1}{1\!+\!c}} \\ \vphantom{\dfrac{-1}{1\!+\!c}} & \lambda \\ & & \vphantom{c^2}\lambda \end{matrix}\\
\end{array} \right|, \nonumber
\end{align}
which can be reduced to a 3x3 matrix using the Schur's complement:
\begin{align}
& |\lambda 1_{6\times 6} - A''| = \begin{vmatrix} A_{11}'' & A_{12}'' \\ A_{21}'' & A_{22}'' \end{vmatrix} = \begin{vmatrix} A_{11}'' & A_{12}'' \\ & A_{22}'' - A_{21}''A_{11}''^{-1}A_{12}'' \end{vmatrix} \nonumber\\
& = |A_{11}''(A_{22}'' - A_{21}''A_{11}''^{-1}A_{12}'')|. \label{eqn:detApp_1}
\end{align}
The determinant of $A_{11}''$ is simply $\lambda^3$, while $-A_{21}''A_{11}''^{-1} A_{21}''$ reduces to
\begin{align}
& \begin{bmatrix} \dfrac{1}{1\!+\!c} \\ & \dfrac{1}{1\!+\!c} \\ 1 & s & c\!-\!c^2 \end{bmatrix} \frac{1_{3\times 3}}{\lambda} \begin{bmatrix} 1\!-\!c & s & 1 \\ -s & 1\!-\!c & s/c \\ & \dfrac{-2s}{1\!+\!c} & \dfrac{-1}{c(1\!+\!c)} \end{bmatrix} \nonumber\\
& = \frac{1}{\lambda} \begin{bmatrix} \dfrac{1}{1\!+\!c} \\ & \dfrac{1}{1\!+\!c} \\ 1 & s & \dfrac{c\!-\!c^2}{1\!+\!c} \end{bmatrix} \begin{bmatrix} 1\!-\!c & s & 1 \\ -s & 1\!-\!c & s/c \\ & -2s & -1/c \end{bmatrix} \nonumber\\
& = \frac{1/\lambda}{1\!+\!c} \begin{bmatrix} 1 \\ & 1 \\ 1\!+\!c & s(1\!+\!c) & c\!-\!c^2 \end{bmatrix} \begin{bmatrix} 1\!-\!c & s & 1 \\ -s & 1\!-\!c & s/c \\ & -2s & -1/c \end{bmatrix} \nonumber\\
& = \frac{1/\lambda}{1\!+\!c} \begin{bmatrix} 1\!-\!c & s & 1 \\ -s & 1\!-\!c & s/c \\ c(c^2\!-\!1) & s(2\!-\!c\!+\!c^2) & \dfrac{1}{c}\!+\!1\!+\!c\!-\!c^2 \end{bmatrix}, \nonumber
\end{align}
yielding the characteristic polynomial
\begin{align}
&|A_{11}''(A_{22}'' - A_{21}''A_{11}''^{-1}A_{12}'')| \nonumber\\
& = \left|\lambda^2 1_{3\times 3} + \frac{1}{1\!+\!c} \begin{bmatrix} 1\!-\!c & s & 1 \\ -s & 1\!-\!c & s/c \\ c^3\!-\!c & s(2\!-\!c\!+\!c^2) & \dfrac{1}{c}\!+\!1\!+\!c\!-\!c^2 \end{bmatrix}\right| \nonumber\\
 & = \begin{vmatrix} \lambda^2 + \dfrac{1\!-\!c}{1\!+\!c} & \dfrac{s}{1\!+\!c} & \dfrac{1}{1\!+\!c} \\ \dfrac{-s}{1\!+\!c} & \lambda^2 + \dfrac{1\!-\!c}{1\!+\!c} & \dfrac{s}{c(1\!+\!c)} \\ c^2\!-\!c & \dfrac{s(2\!-\!c\!+\!c^2)}{1\!+\!c} & \lambda^2 + \dfrac{1\!+\!c\!+\!c^2\!-\!c^3}{c(1\!+\!c)} \end{vmatrix} \nonumber\\
 & = \lambda^6 + \lambda^4\, \dfrac{1+3c-c^2-c^3}{c(1+c)} + \lambda^2\, \dfrac{(1-c)(3+3c-2c^2)}{(1+c)^2} \nonumber\\
 & \qquad + \dfrac{2c(1-c)^2}{(1+c)^3}. \label{eqn:cpoly}
\end{align}
The coefficients of $\lambda^6$ and $\lambda^4$ are straightforward to obtain. To derive the constant term $|A''|$, we may multiply the second row in $A_{11}''(A_{22}'' - A_{21}''A_{11}''^{-1}A_{21}'')$ by $s$ and divide the second column by $s$ to obtain a common factor of $1-c$ in the first column, and also accumulate the common factors of $1+c$ in all the denominators, to obtain
\begin{align}
 & |A''| = \begin{vmatrix} \dfrac{1\!-\!c}{1\!+\!c} & \dfrac{s}{1\!+\!c} & \dfrac{1}{1\!+\!c} \\ \dfrac{-s}{1\!+\!c} & \dfrac{1\!-\!c}{1\!+\!c} & \dfrac{s}{c(1\!+\!c)} \\ c^2\!-\!c & \dfrac{s(2\!-\!c\!+\!c^2)}{1\!+\!c} & \dfrac{1\!+\!c\!+\!c^2\!-\!c^3}{c(1\!+\!c)} \end{vmatrix} \nonumber\\
& = \frac{1}{(1+c)^3} \begin{vmatrix} 1\!-\!c & s & 1 \\ -s & 1\!-\!c & s/c \\ c^3\!-\!c & s(2\!-\!c\!+\!c^2) & 1/c\!+\!1\!+\!c\!-\!c^2 \end{vmatrix} \nonumber\\
& = \frac{1}{(1+c)^3} \begin{vmatrix} 1\!-\!c & 1 & 1 \\ -s^2 & 1\!-\!c & s^2/c \\ c^3\!-\!c & 2\!-\!c\!+\!c^2 & 1/c\!+\!1\!+\!c\!-\!c^2 \end{vmatrix} \nonumber\\
& = \frac{1-c}{(1+c)^3} \begin{vmatrix} 1 & 1 & 1 \\ -1\!-\!c & 1\!-\!c & 1/c\!-\!c \\ -c\!-\!c^2 & 2\!-\!c\!+\!c^2 & 1/c\!+\!1\!+\!c\!-\!c^2 \end{vmatrix} \nonumber\\
& = \frac{1-c}{(1+c)^3} \begin{vmatrix} 1 \\ -1\!-\!c & 2 & 1/c\!+\!1 \\ -c\!-\!c^2 & 2\!\!+\!2c^2 & 1/c\!+\!1\!+\!2c\! \end{vmatrix} = \frac{2c(1-c)^2}{(1+c)^3}. \nonumber
\end{align}
The coefficient of $\lambda^2$ is the trace of the adjugate matrix:
\begin{align}
& \frac{1}{(1+c)^2} \left[2(1-c)(1/c+1+c-c^2) \right. \nonumber\\
& \qquad \left. - (1/c-c)(2-c+c^2) - (c^3-c) + 2(1-c)\right] \nonumber\\
 & \quad = \frac{(1 - c)(3 + 3c - 2c^2)}{(1+c)^2} = \frac{3-5c^2+2c^3}{(1+c)^2}. \nonumber
\end{align}

\end{document}